\documentclass{ifacconf}

\newtheorem{theorem}{Theorem}

\newtheorem{problem}{Problem}

\usepackage{graphicx}      
\usepackage{natbib}        
\usepackage{amsmath}
\usepackage{algorithm}
\usepackage{algorithmic}
\usepackage{amssymb}
\usepackage{amsmath}
\usepackage{amsfonts}
\usepackage{amssymb}

\begin{document}
\begin{frontmatter}

\title{Co-design of controller and routing redundancy over a wireless network\thanksref{footnoteinfo}} 

\thanks[footnoteinfo]{The research leading to these results has received funding from the Italian Government under Cipe resolution n.135 (Dec. 21, 2012), project \emph{INnovating City Planning through Information and Communication Technologies} (INCIPICT).}

\author{G.D. Di Girolamo, A. D'Innocenzo and M.D. Di Benedetto}

\address{Department of Information Engineering, Computer Science and Mathematics, Center of Excellence DEWS, Univ. of L'Aquila, Italy}

\begin{abstract}                
In this paper we investigate the exploitation of redundancy when routing actuation data to a discrete-time LTI system connected to the controller via a wireless network affected by packet drops. We assume that actuation packets can be delivered from the controller to the actuator via multiple paths, each associated with a delay and a packet loss probability. We show that the joint design of controller gain and routing redundancy exploitation can tremendously improve the control performance. To achieve this goal we set up and solve a LQR problem for a class of systems that extends discrete-time Markov Jump Linear Systems, in that both continuous and discrete control signals can be actuated.
\end{abstract}


\end{frontmatter}


\section{Introduction}\label{secIntro}

\begin{figure}[ht]
\begin{center}
\vspace{-0.5cm}
\includegraphics[width=0.5\textwidth]{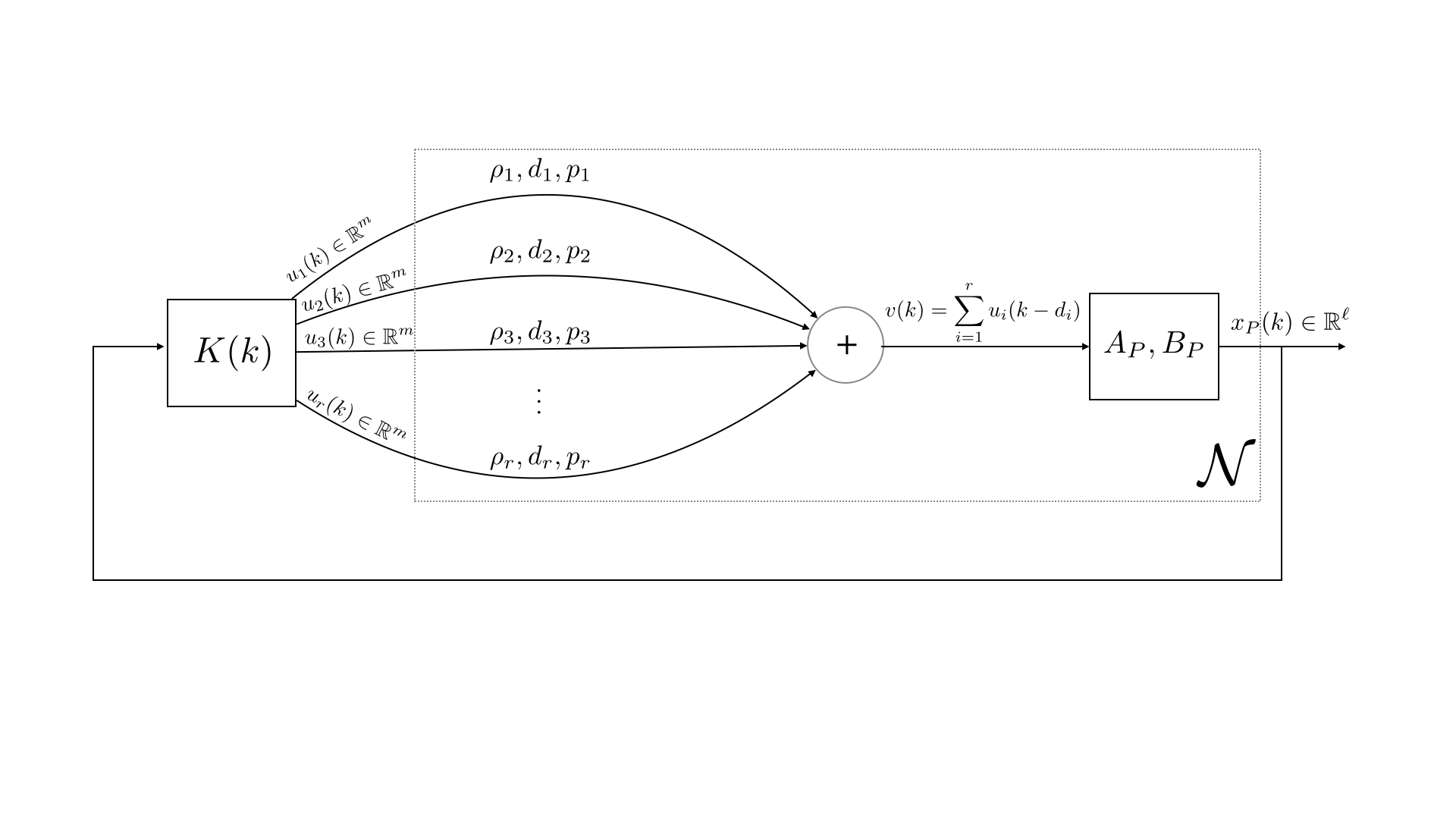}
\vspace{-2cm}
\caption{State feedback control scheme.}\label{figMCNScheme}
\end{center}
\end{figure}

Wireless control networks (WCN) are distributed control systems where the communication between sensors, actuators, and computational units is supported by a wireless communication network. The use of WCN in industrial automation results in flexible architectures and generally reduces installation, debugging, diagnostic and maintenance costs with respect to wired networks (see e.g. \cite{akyildiz_wireless_2004}, \cite{SongIECON2010} and references therein). However modeling, analysis and design of (wireless) networked control systems are challenging open research problems since they require to take into account the joint dynamics of physical systems, communication protocols and network infrastructures. Recently, a huge effort has been made in scientific research on NCSs, see e.g. \cite{Astrom97j1},~\cite{Zhang2001},~\cite{walsh_stability_2002},~\cite{Arzen06},~\cite{TabbaraTAC2007},~\cite{Hespanha2007},~\cite{MurrayTAC2009},~\cite{HeemelsTAC11},~\cite{AlurTAC11},~\cite{PajicTAC2011},~\cite{DInnocenzoTAC13} and references therein for a general overview. 

To make a WCN robust to packets losses redundancy in data routing can be used. One approach to exploit this redundancy is relaying data via multiple paths and then appropriately combining them, which is reminiscent of network coding. In \cite{SmarraECC15} we considered a state-feedback control loop as in Figure \ref{figMCNScheme} where multiple copies $u_1(k)=u_2(k)=\ldots=u_r(k)=Kx_P(k)$ of the same actuation data are sent from the controller to the plant via $r$ different routing paths $\{\rho_i\}_{i=1}^r$ each characterised by a delay $d_i$ and a packet losses probability $p_i$. We assumed that the time-invariant controller gain $K$ is designed via classical methods to assign the eigenvalues of the closed-loop system in the nominal case, i.e. when the effect of packet losses is not considered. We also assumed that the actuator computes a linear combination $\sum_{i=1}^{r}\gamma_i u_i(k-d_i)$ of the data incoming from different routing paths, and we provided a suboptimal algorithm to compute the optimal weights $\gamma_i$ that maximize a metric induced by the notion of Mean Square Stability. In this paper we continue the research line started in \cite{SmarraECC15} and provide novel results that strongly improve the controller performance. The first difference is motivated by the following consideration: when routing redundancy is exploited in communication systems the objective is to relay some information, and thus we send to the network the same packet $u_1(k)=\ldots=u_r(k)$ and try to extract from the corrupted received packets the original information; in our case the objective is to increase the control performance, as a consequence the actuation packets $u_1(k), \ldots, u_r(k)$ must not be necessarily equal. In this paper we perform \emph{controller and routing redundancy co-design} by designing the time-varying matrix $\mathbb{R}^{rm \times n} \ni K(k) \doteq [K_1(k) \in \mathbb{R}^{m \times n}; \cdots; K_r(k)\in \mathbb{R}^{m \times n}]$. Note that the problem formulation in \cite{SmarraECC15} is a special case of the above definition where $K(k) = [\gamma_1 K^*; \cdots; \gamma_r K^*]$ and $K^*$ is designed for the nominal case. As a further improvement, while in \cite{SmarraECC15} we optimise a metric based on the notion of Mean Square Stability (i.e. only taking into account the steady state behavior) we consider here a more complex control specification that also takes into account the transient behavior by setting up a finite-horizon LQR problem. In \cite{MHNAUT12} the authors also consider redundant data transmission over a set of paths characterised by i.i.d. Bernoulli probabilities of packet losses, but for a more restricted scenario with respect to ours because they assume that all paths are associated with the same delay, the packet loss events are measurable, the controller is designed for the nominal case (i.e. without considering the effect of packet losses) and redundant data combination is not modeled/designed. Their focus is on deciding how many redundant copies of a packet should be transmitted at each sampling time and what benefits can be drawn from this: besides the fact that our model is more general, we also address the more general problem of co-designing controller and routing redundancy.

In Section \ref{secModel} we define a network modeling framework that allows co-design of controller and routing redundancy while taking into account the effect of packet losses and in Section \ref{secProblem} we provide our LQR problem formulation. 

In Section \ref{secStaticRouting} we assume that the set of paths used at each time instant to send actuation data to the actuator has be designed a priori: we will call this approach static routing redundancy. We setup the problem of co-desiging the controller gain and the routing redundancy parameters as a LQR problem for a Discrete-Time Markov-Jump Linear System (dtMJLS) where the discrete mode (which correpsond to the occurrence of packet losses) is unmeasurable and evolves according to a sequence of i.i.d. random variables. The latter assumption, widely adopted for several communication systems, makes particularly sense in our framework since exploitation of redundant data is well known to be very effective especially in the case when the reduntant paths used for data relay are uncorrelated from the point of view of the communication channels' characteristics. The proof for solving the LQR problem for such a system, being the discrete state unmeasurable, is an extension of the solution for dtMJLS in \cite{costa_discrete-time_2005} and can be derived without much difficulty thanks to the i.i.d. assumption. Note that a similar problem has been addressed in \cite{BarasProcIFAC}, with the assumption that the discrete state is measurable with a one step delay, and solved without proof: we believe that the proof in this case is very close to ours, but since we were unable to find it in the scientific literature we provide one in this paper for the sake of completeness. We apply our optimal solution to a simple example where actuation data can be sent to the actuator via two paths: the first characterised by short delay (i.e. fast reaction to perturbations) and high probability of packet losses (i.e. low reliablity); the second characterised by long delay (i.e. slow reaction to perturbations) and $0$ probability of packet losses (i.e. perfect reliablity). Note that such situation can often occur in realistic cases: one example is a multi-hop wireless network where we can reach the destination via a single long hop (short delay, high packet loss probability) or via a path of very short multiple hops (high delay, low packet loss probability); another example is a service provider network where we can reach the destination via the shortest yet congested path of routers (short delay, high packet loss probability) or via a longer uncongested path of routers (high delay, low packet loss probability). In the above situations, using only the first path is clearly not a good idea since the closed loop system may easily become unstable. Using only the second path is the optimal solution to maximise bandwidht, i.e. optimal from the point of view of communication theory: however, due to the high delay, the control system is not reactive to perturbations. The main idea that motivates this paper is based on the intuition that we could use both paths simultaneously, exploiting the fast reaction of the first path and the high reliability of the second path in an optimal way taking into account the plant dynamics. Our Monte Carlo (MC) simultations show that routing actuation data on both paths simultaneously and applying our optimal solution, we can tremendously improve the performance of both single-path solutions from the point of view of control performance.

In Section \ref{secDynamicRouting} we consider a much more complicated problem: we assume that the set of paths used at each time instant to send data to the actuator can be controlled, i.e. the choice of redundant routing paths is also a control variable: we will call this approach dynamic routing redundancy. We setup the problem of co-desiging the controller gain, routing redundancy parameters and paths as a LQR problem for a class of systems that includes dtMJLS as a special case and provide, as the main theoretical contribution of this paper, a recursive solution that is optimal for a certain set of initial conditions, which we define in closed form. More precisely, our model is a dtMJLS where we can also apply a discrete control that, choosing some of the system matrices, models the choice of the routing at each time step. A similar model has been considered in \cite{VargasCDC2010}, where only a sub-optimal solution is provided based on a conservative approximation (see the proof of this paper for more details). Similar problems have also been considered, but for determinstic models, in \cite{AbateACC2009} and \cite{BemporadAutomatica2005}. Finally, in \cite{PappasTAC2014} a different LQR optmisation problem is considered where the discrete and continuous control signals are independent and can be designed separately, which is not the case in our problem.

%
%


\section{Network Modeling}\label{secModel}

Consider a state-feedback networked control loop as in Figure \ref{figMCNScheme} where the communication between the controller and the actuator can be performed via a set of $r$ routing paths $\{\rho_i\}_{i=1}^r$ in a wireless multi-hop communication network. Each path $\rho_i$ is characterised by a delay $d_i \in \mathbb{N^+}$ and a packet loss probability $p_i \in [0,1]$ that represents the probability that the packet transmitted on that path will not reach the actuator due to communication failure. In particular, we define for each path $\rho_i$ the stochastic process $\sigma_i(k) \in \{0,1\}$, with $\sigma_i(k)=0$ if the packet expected to arrive via the routing path $\rho_i$ at time $k$ suffered a packet drop and $\sigma_i(k)=1$ if the packet is succesfully received at time $k$. We assume that $\sigma_i(k)$ is a sequence of i.i.d. random variables, each characterised by a Bernoulli distribution with probability measure $\mathbb{P}[\sigma_i(k)=0]=p_i$. We also assume here that the events of occurrence of packet losses in the different paths are i.i.d.: as a consequence the stochastic process $\sigma(k) \doteq [\sigma_1(k), \ldots, \sigma_r(k)]'$ is a vector of i.i.d. random variables, where $\sigma(k)$ can assume $2^r$ values. We also assume that the controller cannot measure the signal $\sigma(k)$, i.e. it is not possible to measure the occurrence of packet losses. The case when the occurrence of packet losses is measurable with finite delay thanks to acknowledgement packets will be investigated in future work. We assume that, in general, the controller can decide for each time instant $k$ the set of paths where data will be sent: i.e., the controller can decide to send data at time $k$ on all paths, on a subset of paths, on one path, or even not to send any data. To this aim we define for each path $i$ the discrete control signal $a_i(k) \in \{0,1\}$, with $a_i(k)=1$ if the controller decides to send a packet via the routing path $i$ at time $k$, and $a_i(k)=0$ if no packet is sent via path $i$ at time $k$. We define the discrete control signal $a(k) \doteq [a_1(k), \ldots, a_r(k)]'$, where $a(k)$ can be choosen among $2^r$ different values.

Let the plant be a discrete-time LTI system described by the matrices $A_P \in \mathbb{R}^{\ell \times \ell}$ , $B_P \in \mathbb{R}^{\ell \times m}$, we define the dynamics of the networked system as follows:
\begin{equation}\label{eqMainModelNetwork}
   \begin{cases}
   x(k+1) = A_{\sigma(k)}x(k) + B_{a(k)}u(k)\\
   y(k) = x(k)
   \end{cases}
\end{equation}
with
\tiny
\begin{align*}
A_{\sigma(k)} &= \left[\begin{matrix}
  A_P & \Lambda_{1}(\sigma(k)) & \Lambda_{2}(\sigma(k)) & \cdots & \Lambda_{r}(\sigma(k)) \\
  0 & \Gamma_{1} & 0 & \cdots & 0 \\
  0 & 0 & \Gamma_{2} & \cdots & 0 \\
  \vdots & \vdots & \vdots & \ddots & \vdots \\
  0 & 0 & 0 & \cdots & \Gamma_{r}
\end{matrix}\right] \in \mathbb R^{\ell+\nu(r) \times \ell+\nu(r)},\\
B_{a(k)} &= \left[\begin{matrix}
0 & 0 & \cdots & 0\\
a_1(k) I_m \otimes \textbf{e}_{d_1} & 0 & \cdots & 0 \\
0 & a_2(k) I_m \otimes \textbf{e}_{d_2} & \cdots & 0 \\
\vdots & \vdots & \ddots & \vdots \\
0 & 0 & \cdots & a_r(k) I_m \otimes \textbf{e}_{d_r} \\
\end{matrix}\right] \in \mathbb R^{\ell + \nu(r) \times mr},
\end{align*}
\normalsize
with
\tiny
\begin{align*}
&\Lambda_i(\sigma(k)) \doteq
\sigma_i(k) \left[\begin{matrix}
  B_P & 0 & \cdots & 0
\end{matrix}\right] \in \mathbb R^{\ell \times md_i},\\
&\Gamma_{i} \doteq \left[\begin{matrix}
  0 & I_m & \cdots & 0 & 0 \\
  \vdots & \vdots & \ddots & \vdots & \vdots \\
  0 & 0 & \cdots & I_m & 0 \\
  0 & 0 & \cdots & 0 & I_m \\
  0 & 0 & \cdots & 0 & 0
\end{matrix}\right] \in \mathbb R^{{md_i} \times {md_i}},\\
\end{align*}
\normalsize
and with $\nu(i) \doteq m\sum\limits_{j=1}^{i}d_j$, $I_m$ the $m$-dimensional identity matrix, $\textbf{e}_{i}$ a column vector of appropriate dimension with all zero entries except the $i-th$ entry equal to 1, and $\otimes$ the Kronecker product. System \eqref{eqMainModelNetwork} is more general than dtMJLSs: in particular, when $\forall k \geq 0, a(k) = a_k$ (i.e. the routing is designed a priori) System \eqref{eqMainModelNetwork} is a dtMJLSs. Note that in our feedback scheme we assume that the controller can measure the whole state $x(k) = [x_P(k)' x_N(k)']'$ of \eqref{eqMainModelNetwork}, where $x_P(k) \in \mathbb{R}^\ell$ is the state of the plant and $x_N(k) \in \mathbb{R}^{\nu(r)}$ are state variables modeling the flow of packets via all routing paths. We assume that the controller can measure the state $x_P(k)$ of the plant via sensors and defer to future work the output-feedback case. Also, the controller is aware of the current and past actuation signals $u(k)$ that have been sent to the actuator, as well as of the current and past signals $a(k)$ (which is either constant for any $k$ or choosen by the controller): as a consequence the controller has direct access to the state of $x_N(k)$, which models the actuation commands that are expected to arrive at the actuator, but is not aware of their arrival to the actuator since $\sigma(k)$ is not measurable. 


\section{Problem Formulation}\label{secProblem}

The network modeling framework in the previous section is a special case of the following mathematical framework, which is the one we will use in the rest of the paper:
\begin{equation}\label{eqMainModel}
   \begin{cases}
   x(k+1) = A_{\sigma(k)}x(k) + B_{a(k)}u(k)\\
   y(k) = x(k)
   \end{cases},
\end{equation}
where $x(k) \in \mathbb{R}^n, u(k) \in \mathbb{R}^{mr}, \sigma(k) \in \{1, \ldots, q\} \doteq \Sigma, a(k) \in \{1, \ldots, p\} \doteq A$, and $\sigma(k), k \geq 0$ is a sequence of i.i.d. random variables such that $\mathbb{P}[\sigma(k) = i] = \pi_i$ for any $i \in \Sigma, k \geq 0$.

\begin{problem}\label{probMain}
Given System \eqref{eqMainModel}, design for any $k \in \{0,\ldots,N-1\}$ an optimal state-feedback control policy $a^*(x(k)), u^*(k) = K^*(x(k))x(k)$, with $a^*(x(k)): \mathbb{R}^n \rightarrow A$ and $K^*(x(k))$ a $m \times n$ matrix of reals, minimizing the following objective function:
\tiny
\begin{align*}
&J(x(0), u(0), a(0))=E \left\{ \sum_{k=0}^{N-1} \left(x'(k)Mx(k)+u'(k)Ru(k)\right) + x'(N)Qx(N) \big| \aleph_{0}\right\}
\end{align*}
\normalsize
where $\aleph_{k}$ is the sigma algebra generated by $x(0),...,x(k)$. 
\end{problem}


\section{Co-design of controller and static routing redundancy}\label{secStaticRouting}

\begin{figure}[h!]
\begin{center}
\includegraphics[width=0.5\textwidth]{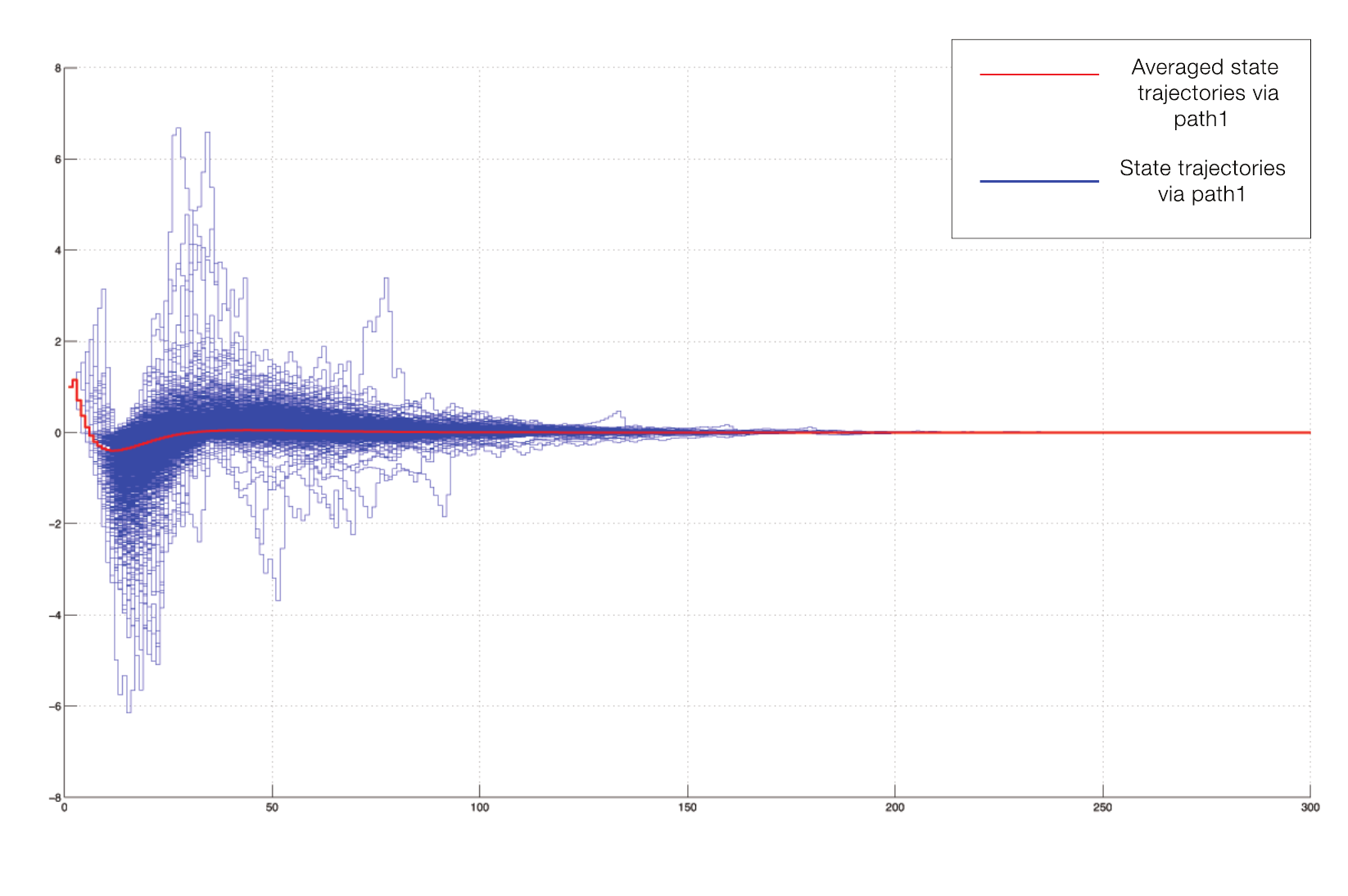}
\caption{State trajectories routing only via path $\rho_1$ (blue) and their average (red).}\label{figPath1}
\end{center}
\end{figure}

\begin{figure}[h!]
\begin{center}
\includegraphics[width=0.5\textwidth]{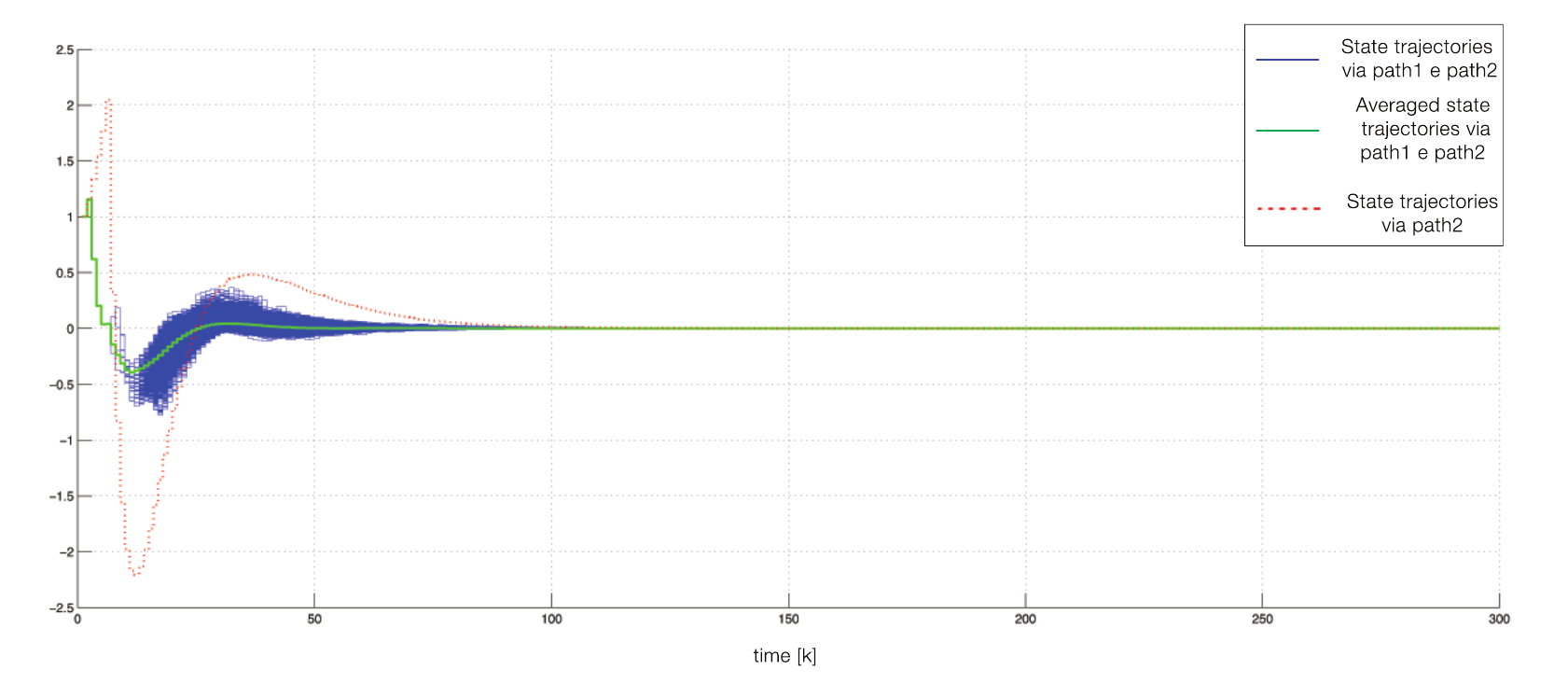}
\caption{State trajectory routing only via path $\rho_2$ (red dashed); state trajectories routing via both paths $\rho_1$ and $\rho_2$ (blue) and their average (green).}\label{figPaths12}
\end{center}
\end{figure}

\begin{figure}[h!]
\begin{center}
\includegraphics[width=0.5\textwidth]{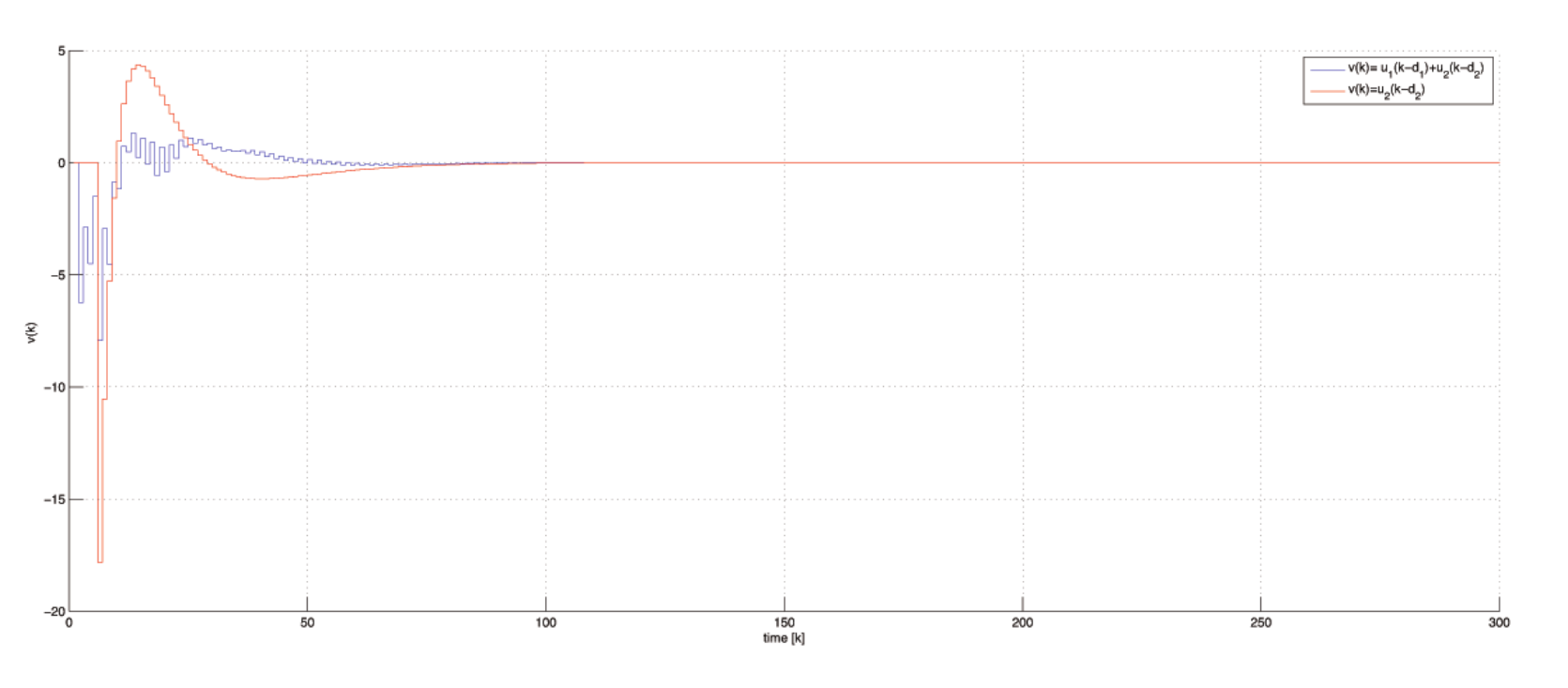}
\caption{Actuation signal $v_2(k) = u_2(k-5)$ (red) routing via path $\rho_2$; averaged actuation signal $v_{1,2}(k) = u_1(k-1) + u_2(k-5)$ (blue) routing via both paths.}\label{figControlSignals}
\end{center}
\end{figure}

In this Section we address Problem \ref{probMain} assuming that the routing has beed designed a priori.
\begin{theorem}\label{thStaticRoutingResult}
Given System \eqref{eqMainModel} and a routing policy defined by $\forall k \geq 0, a(k) = a_k \in A$, the optimal solution of Problem \ref{probMain} is given by a sequence $K^*(k)$ with $k =0,\ldots,N-1$.

\emph{Proof:} The proof is constructive and shows how to compute $K^*(k)$ for any $k =0,\ldots,N-1$.We start from the classical Belmann optimization formulation (\cite{LancasterRiccati}, \cite{BertsekasOptStocCon}):
\tiny
\begin{equation}
   \begin{cases}
   J(x(k), u(k)) =\\
   \min\limits_{a(k), u(k)} E \{x'(k)Mx(k)+u'(k)Ru(k)+J(x(k+1), u(k+1)) | \aleph_{k} \}\\
   \\
   J(x(N) , u(N)) = J(x(N)) = x'(N) E\{Q | \aleph_{N}\} x(N) = x'(N) P(N) x(N)
   \end{cases}
\end{equation}
\normalsize
where $P(N) \doteq Q$ is a symmetric matrix and $J(x(k), u(k))$ is the cost-to-go function at time $k$. Let us write the cost-to-go function at step $N-1$:
\footnotesize
\begin{align*}
& J(x(N-1),u(N-1))=\\
& \min_{u(N-1)} E \Big\{x'(N-1)Mx(N-1)+u'(N-1)Ru(N-1) +\\
& + x'(N)P(N)x(N) \Big| \aleph_{N-1} \Big\}=\\
& \min_{u(N-1)} E \Big\{x'(N-1)Mx(N-1)+u'(N-1)Ru(N-1) +\\
& (x'(N-1)A'_{\sigma(N-1)}+u'(N-1)B_{a_{N-1}}')P(N)(A_{\sigma(N-1)}x(N-1)+\\
& B_{a_{N-1}}u(N-1)) \Big| \aleph_{N-1} \Big\}
\end{align*}
\normalsize
By linearity of the expected value we can write:
\tiny
\begin{align}\label{LQRsplitStaticRouting}
&J(x(N-1),u(N-1))=\notag\\
&\min\limits_{u(N-1)} \Big\{ E \{x'(N-1)Mx(N-1) | \aleph_{N-1} \}+\notag\\
&E \{u'(N-1)Ru(N-1)   | \aleph_{N-1}\} +\notag\\
&E \{(x'(N-1)A'_{\sigma(N-1)}P(N)A_{\sigma(N-1)}x(N-1) | \aleph_{N-1}\} +\notag\\
&E \{(x'(N-1)A'_{\sigma(N-1)}P(N)B_{a_{N-1}}u(N-1) | \aleph_{N-1}\} +\notag\\
&E \{(u'(N-1)B_{a_{N-1}}'P(N)A_{\sigma(N-1)}x(N-1)  | \aleph_{N-1}\} +\notag\\
&E \{u'(N-1))B_{a_{N-1}}'P(N)B_{a_{N-1}}u(N-1) |\aleph_{N-1}\} \Big\}
\end{align}
\normalsize
Let us consider each addend of the rightside of \eqref{LQRsplitStaticRouting}. Since $M$ and $R$ are constant matrices, $x(N-1)$ is $\aleph_{N-1}$-measurable, and $u(N-1)$ is not a random variable since it is the input that we choose to apply to the system at time $N-1$, the first two addends can be written as
\tiny
\begin{align*}
&E \{x'(N-1)Mx(N-1) +u'(N-1)Ru(N-1)| \aleph_{N-1}\}= \\
&x'(N-1)Mx(N-1) +u'(N-1)Ru(N-1).
\end{align*}
\normalsize
Moreover, $x'(N-1)Mx(N-1)$ does not depend on $u(N-1)$ and can be moved out of the min operator. The third addend can be written as
\begin{align*}
&E \{(x'(N-1)A'_{\sigma(N-1)}P(N)A_{\sigma(N-1)}x(N-1) | \aleph_{N-1}\}=\\
&x'(N-1) E \{A'_{\sigma(N-1)}P(N)A_{\sigma(N-1)}| \aleph_{N-1}\} x(N-1) \doteq\\
&x'(N-1) \Phi(N-1)(N-1),
\end{align*} 
\normalsize
where $\Phi(N-1) = \sum\limits_{i=1}^{q} A'_iP(N)A_i \pi_i$. Since $x(N-1)$ is $\aleph_{N-1}$-measurable, $P(N)$ is symmetric and $\sigma(k)$ are i.i.d., the sum of the fourth and fifth addends can be written as
\begin{align*}
&2E \{(x'(N-1)A'_{\sigma(N-1)}P(N)B_{a_{N-1}}u(N-1) | \aleph_{N-1}\}=\\
&2x'(N-1)E \{ A'_{\sigma(N-1)}P(N)| \aleph_{N-1}\}B_{a_{N-1}}u(N-1) =\\
&2x'(N-1)E \{ A'_{\sigma(N-1)} \} P(N) B_{a_{N-1}}u(N-1) =\\
&2x'(N-1) \bar A' P(N) B_{a_{N-1}}u(N-1),
\end{align*}
where $\bar A \doteq E\{A_{\sigma(N-1)}\} = \sum_{i=1}^q A_i \pi_i$. The last addend can be written as:
\begin{align*}
&E \{u'(N-1))B_{a_{N-1}}'P(N)B_{a_{N-1}}u(N-1) |\aleph_{N-1}\}=\\
&u'(N-1))B_{a_{N-1}}' P(N) B_{a_{N-1}}u(N-1).
\end{align*}
We can now rewrite \eqref{LQRsplitStaticRouting} as follows:
\tiny
\begin{align} \label{costo_n-1StaticRouting}
&J(x(N-1),u(N-1))= x'(N-1)\left[M+ \Phi(N-1) \right]x(N-1)+\notag\\
&\min\limits_{u(N-1)} \Big\{ 2x'(N-1) \bar A' P(N) B_{a_{N-1}}u(N-1)+\notag\\ 
& u'(N-1))B_{a_{N-1}}' P(N) B_{a_{N-1}}u(N-1)+u'(N-1)Ru(N-1)\Big\}
\end{align}
\normalsize
We can compute the minimun of the above function with respect to $u(N-1)$ by equaling the derivative with respect to $u(N-1)$ to $0$:
\tiny
\begin{equation}
	\begin{split}
		2x'(N-1) \bar A' P(N) B_{a_{N-1}}	+ 2u'(N-1)\left[ B_{a_{N-1}}' P(N) B_{a_{N-1}}+R \right] =0,
	\end{split}
\end{equation}
\normalsize
which gets to
\tiny
\begin{align*}
&u(N-1)= -\left[R+ B_{a_{N-1}}' P(N) B_{a_{N-1}} \right]^{-1} B_{a_{N-1}}' P(N) \bar A x(N-1).
\end{align*}
\normalsize
We have thus obtained a linear feedback of the state given by
\begin{equation}
u(N-1)=K^*(N-1)x(N-1) 
\end{equation}
By replacing the expression of $u(N-1)$ in the cost function \eqref{costo_n-1StaticRouting} it is possible to obtain
\tiny
\begin{align*}
& J(x(N-1)) =\\
& x'(N-1)\Big[M+\Phi(N-1) + \bar A'  P(N) B_{a_{N-1}} \left(R+B_{a_{N-1}}'P(N) B_{a_{N-1}} \right)^{-1} \cdot \\
& B_{a_{N-1}}' P(N) \bar A \Big] x(N-1) \doteq x'(N-1) P(N-1) x(N-1).
\end{align*}
\normalsize 
By iterating this procedure for a generic $k$ it is possible to obtain the expression for $K^*(k)$
\begin{equation} \label{controllore}
K^*(k) =	-\left[R+ B_{a_{k}}' P(k+1) B_{a_{k}} \right]^{-1} B_{a_{k}}'P(k+1) \bar A.
\end{equation}
This concludes the proof.$\qed$
\end{theorem}
We apply the above theorem to the following System \eqref{eqMainModelNetwork} characterised by a 4-dimensional unstable randomly generated plant
\tiny
\begin{align*}
A_P &= \left[\begin{matrix}
    1.1062 &  -1.0535  &  0.7944  & -0.4543\\
    0.0202  & -0.0654  &  0.9697  & -0.6888\\
    0.1131 &  -0.5755  &  1.7434 &  -0.7174\\
    0.0745 &  -0.2565  &  0.2999  &  0.7252\\
\end{matrix}\right],
B_P = \left[\begin{matrix}
   -0.1880\\
    0.0182\\
    0.1223\\
    0.2066\\
\end{matrix}\right],
\end{align*}
\normalsize
and by a wireless network carachterized by two paths: $\rho_1$ with packet loss probability $p_1=0.25$ and delay $d_1=1$ and $\rho_2$ with packet loss probability $p_2 = 0$ and delay $d_2=5$. We applied Theorem \ref{thStaticRoutingResult} to three simple routing strategies, i.e. using for all time instants only path $\rho_1$, only path $\rho_2$ or both paths simultaneously. We computed our solution for a time horizon $N=300$ and defined the weight matrices $M,R,Q$ and initial condition $x(0)$ as follows:
\tiny
\begin{align*}
&\text{If routing via both } \rho_1,\rho_2:\\
&M=Q = \left[\begin{matrix}
  I_4 & 0\\
  0       & 0
\end{matrix}\right] \in \mathbb R^{10 \times 10}, R = I_2, x(0) = (1,1,1,1,0,\ldots,0) \in \mathbb R^{10};\\
&\text{if routing via } \rho_1:\\
&M=Q= \left[\begin{matrix}
  I_4 & 0\\
  0       & 0
\end{matrix}\right] \in \mathbb R^{5 \times 5},\quad R = 1, \quad x(0) = (1,1,1,1,0,\ldots,0) \in \mathbb R^{5};\\
&\text{if routing via } \rho_2:\\
&M=Q = \left[\begin{matrix}
  I_4 & 0\\
  0       & 0
\end{matrix}\right] \in \mathbb R^{9 \times 9},\quad R = 1, \quad x(0) = (1,1,1,1,0,\ldots,0) \in \mathbb R^{9}.
\end{align*}
\normalsize
For each routing strategy we performed $5K$ MC simulations of the state trajectories. Figure \ref{figPath1} shows the trajectories of the first component of the extended state vector when only path $\rho_1$ is used. The system can be stabilized, but the variance of the trajectories is high. In many other simulations the whole system can't be stabilized only using path $\rho_1$, many state trajectories are unstable because of the very high packet loss probability. This routing policy is clearly a bad choice.
Figure \ref{figPaths12} shows the trajectories when only path $\rho_2$ is used (Red) and when both paths $\rho_1$ and $\rho_2$ are used (Blue and green). Routing data only to path $\rho_2$ clearly generates always the same trajectory since $p_2 = 0$. The system trajectory is stable but the associated cost is quite high because of the large delay, as evidenced by the overshoot and the settling time performances. Figure \ref{figPaths12} evidences that routing data via both paths $\rho_1$ and $\rho_2$ the control performance strongly improves: in particular, the trajectory of the system computed by averaging over all MC simulations is characterised by much smaller overshoot and faster settling time. The single trajectories generated routing data via both paths $\rho_1$ and $\rho_2$ clearly have some variance due to the high packet loss probability $p_1$: however, in the $5K$ MC simulations, the performance of any of the single trajectories is much better than the case when only path $\rho_2$ is used.
Table \ref{tabCosts} shows the tremendous improvement of the controller performance obtained by exploiting both paths and co-designing controller gains and static routing redundancy parameters.
\begin{table}[h] \label{tabCosts}
\begin{center} 
\begin{tabular}{|c|c|}
\hline
& Averaged cost\\ \hline
Via path $\rho_1$ & $\sim 900$\\ \hline
Via path $\rho_2$ & $\sim 250$\\ \hline
Via paths $\rho_1, \rho_2$ & $\sim 100$\\ \hline
\end{tabular}
\end{center}
\caption{Cost averaged over 5K MC sims.}
\end{table}
Finally, the averaged actuation signals $v_2(k) = u_2(k-5)$ routing via path $\rho_2$ and $v_{1,2}(k) = u_1(k-1) + u_2(k-5)$ routing via both paths are shown in Figure \ref{figControlSignals}: note that the actuation cost (i.e. energy) associated to the case when we use only path $\rho_2$ is much larger w.r.t. to the case when we use both paths $\rho_1$ and $\rho_2$.
Using our approach it is possible to compute, for a finite set of predefined routing policies, the associated expected quadratic cost and choose the cheapest policy. To further improve the performance one can dynamically choose, for each packet and according to the plant state measurement, the routing path: this is the problem we address in the next section.

\section{Co-design of controller and dynamic routing redundancy}\label{secDynamicRouting}

In this section we provide a suboptimal solution of Problem \ref{probMain} that is optimal for a certain set of initial conditions, which we define in closed form.
\begin{theorem}\label{thMainResult}
Given System \eqref{eqMainModel}, a solution of Problem \ref{probMain} is given by the sequence, $\forall k =0,\ldots,N-1$, of $a^*(x(k))$, $K^*(x(k))$ defined by a finite partition $\Omega(k) \doteq \{\Omega_{i}(k)\}_{i=1}^{\omega(k)}$ of $\mathbb{R}^n$, each assocated with a pair $a_i, K_i$ such that $a^*(x(k)) = a_i, K^*(x(k)) = K_i \iff x(k) \in \Omega_i(k)$. Any $\Omega_i(k)$ can be defined by a finite set of quadratic inequalities in the form $x'Yx \sim 0$, $\sim \in \{<, \leq, \geq, >\}$ with $Y$ a $n \times n$ symmetric matrix of reals. Moreover, at each step $k$ such solution is optimal for a subset $\iota(k) \subseteq \{0,\ldots,\omega(k)\}$. 

\emph{Proof:} The proof is constructive and shows how to compute each $\Omega_i(k)$. We start from the classical Belmann optimization formulation:
\tiny
\begin{equation}
   \begin{cases}
   J(x(k), u(k), a(k)) =\\
   \min\limits_{a(k), u(k)} E \{x'(k)Mx(k)+u'(k)Ru(k)+J(x(k+1), u(k+1), a(k+1)) | \aleph_{k} \}\\
   \\
   J(x(N) ,u(N), a(N)) = J(x(N)) = x'(N) E\{Q | \aleph_{N}\} x(N) = x'(N) P(N) x(N)
   \end{cases}
\end{equation}
\normalsize
where $P(N) \doteq Q$ is a symmetric matrix and $J(x(k), u(k), a(k))$ is the cost-to-go function at time $k$. We first provide the optimal solution at step $N-1$. Then we provide the optimal solution at step $N-2$ given the optimal solution at step $N-1$. The optimal solution at any other step $k = 1, \ldots, N-3$ can be obtained iterating the same reasoning of step $N-2$.

\textit{Step $N-1$:} Let us write the cost-to-go function at step $N-1$:
\footnotesize
\begin{align*}
& J(x(N-1),u(N-1),a(N-1))=\\
& \min_{a(N-1),u(N-1)} E \Big\{x'(N-1)Mx(N-1)+u'(N-1)Ru(N-1) +\\
& + x'(N)P(N)x(N) \Big| \aleph_{N-1} \Big\}=\\
& \min_{a(N-1),u(N-1)} E \Big\{x'(N-1)Mx(N-1)+u'(N-1)Ru(N-1) +\\
& (x'(N-1)A'_{\sigma(N-1)}+u'(N-1)B_{a(N-1)}')P(N)(A_{\sigma(N-1)}x(N-1)+\\
& B_{a(N-1)}u(N-1)) \Big| \aleph_{N-1} \Big\}
\end{align*}
\normalsize
By linearity of the expected value we can write:
\tiny
\begin{align}\label{LQRsplit}
&J(x(N-1),u(N-1),a(N-1))=\notag\\
&\min\limits_{a(N-1),u(N-1)} \Big\{ E \{x'(N-1)Mx(N-1) | \aleph_{N-1} \}+\notag\\
&E \{u'(N-1)Ru(N-1)   | \aleph_{N-1}\} +\notag\\
&E \{(x'(N-1)A'_{\sigma(N-1)}P(N)A_{\sigma(N-1)}x(N-1) | \aleph_{N-1}\} +\notag\\
&E \{(x'(N-1)A'_{\sigma(N-1)}P(N)B_{a(N-1)}u(N-1) | \aleph_{N-1}\} +\notag\\
&E \{(u'(N-1)B_{a(N-1)}'P(N)A_{\sigma(N-1)}x(N-1)  | \aleph_{N-1}\} +\notag\\
&E \{u'(N-1))B_{a(N-1)}'P(N)B_{a(N-1)}u(N-1) |\aleph_{N-1}\} \Big\}
\end{align}
\normalsize
Let us consider each addend of the rightside of \eqref{LQRsplit}. Since $M$ and $R$ are constant matrices, $x(N-1)$ is $\aleph_{N-1}$-measurable, and $u(N-1)$ is not a random variable since it is the input that we choose to apply to the system at time $N-1$, the first two addends can be written as
\tiny
\begin{align*}
&E \{x'(N-1)Mx(N-1) +u'(N-1)Ru(N-1)| \aleph_{N-1}\}= \\
&x'(N-1)Mx(N-1) +u'(N-1)Ru(N-1).
\end{align*}
\normalsize
Moreover, $x'(N-1)Mx(N-1)$ does not depend on $a(N-1), u(N-1)$ and can be moved out of the min operator. The third addend can be written as
\begin{align*}
&E \{(x'(N-1)A'_{\sigma(N-1)}P(N)A_{\sigma(N-1)}x(N-1) | \aleph_{N-1}\}=\\
&x'(N-1) E \{A'_{\sigma(N-1)}P(N)A_{\sigma(N-1)}| \aleph_{N-1}\} x(N-1) \doteq\\
&x'(N-1) \Phi(N-1)(N-1),
\end{align*} 
\normalsize
where $\Phi(N-1) = \sum\limits_{i=1}^{q} A'_iP(N)A_i \pi_i$. Since $x(N-1)$ is $\aleph_{N-1}$-measurable, $P(N)$ is symmetric and $\sigma(k)$ are i.i.d., the sum of the fourth and fifth addends can be written as
\begin{align*}
&2E \{(x'(N-1)A'_{\sigma(N-1)}P(N)B_{a(N-1)}u(N-1) | \aleph_{N-1}\}=\\
&2x'(N-1)E \{ A'_{\sigma(N-1)}P(N)| \aleph_{N-1}\}B_{a(N-1)}u(N-1) =\\
&2x'(N-1)E \{ A'_{\sigma(N-1)} \} P(N) B_{a(N-1)}u(N-1) =\\
&2x'(N-1) \bar A' P(N) B_{a(N-1)}u(N-1),
\end{align*}
where $\bar A \doteq E\{A_{\sigma(N-1)}\} = \sum_{i=1}^q A_i \pi_i$. The last addend can be written as:
\begin{align*}
&E \{u'(N-1))B_{a(N-1)}'P(N)B_{a(N-1)}u(N-1) |\aleph_{N-1}\}=\\
&u'(N-1))B_{a(N-1)}' P(N) B_{a(N-1)}u(N-1).
\end{align*}
We can now rewrite \eqref{LQRsplit} as follows:
\tiny
\begin{align} \label{costo_n-1}
&J(x(N-1),u(N-1),a(N-1))= x'(N-1)\left[M+ \Phi(N-1) \right]x(N-1)+\notag\\
&\min\limits_{a(N-1),u(N-1)} \Big\{ 2x'(N-1) \bar A' P(N) B_{a(N-1)}u(N-1)+\notag\\ 
& u'(N-1))B_{a(N-1)}' P(N) B_{a(N-1)}u(N-1)+u'(N-1)Ru(N-1)\Big\}
\end{align}
\normalsize
For any given $a(N-1) \in A$ we can compute the minimun of the above function with respect to $u(N-1)$ by equaling the derivative with respect to $u(N-1)$ to $0$:
\tiny
\begin{equation}
	\begin{split}
		2x'(N-1) \bar A' P(N) B_{a(N-1)}	+ 2u'(N-1)\left[ B_{a(N-1)}' P(N) B_{a(N-1)}+R \right] =0,
	\end{split}
\end{equation}
\normalsize
which gets to
\tiny
\begin{align*}
&u_{a(N-1)}(N-1)= -\left[R+ B_{a(N-1)}' P(N) B_{a(N-1)} \right]^{-1} B_{a(N-1)}' P(N) \bar A x(N-1).
\end{align*}
\normalsize
We have thus obtained a linear feedback of the state given by
\begin{equation}
u_{a(N-1)}(N-1)=K_{a(N-1)}(N-1)x(N-1) 
\end{equation}
By replacing the expression of $u_{a(N-1)}(N-1)$ in the cost function \eqref{costo_n-1} it is possible to obtain
\tiny
\begin{align*}
& J(x(N-1), a(N-1)) =\\
& x'(N-1)\Big[M+\Phi(N-1) + \bar A'  P(N) B_{a(N-1)} \left(R+B_{a(N-1)}'P(N) B_{a(N-1)} \right)^{-1} \cdot \\
& B_{a(N-1)}' P(N) \bar A \Big] x(N-1) \doteq x'(N-1) P_{a(N-1)}(N-1) x(N-1).
\end{align*}
\normalsize 
Let us now define a partition of $\mathbb{R}^n$ given by the collection of disjoint sets $\Omega(N-1) \doteq \{\Omega_i(N-1)\}_{i=1}^p$ where each $\Omega_i(N-1)$ is defined by
\tiny
\begin{align*}
\Omega_i(N-1) \doteq \{x(N-1) \in \mathbb{R}^n : i = \arg\min\limits_{a(N-1)} J(x(N-1), a(N-1)) \}.
\end{align*}
\normalsize
$\Omega_i(N-1)$ is the set of all states $x(N-1)$ such that the corresponding optimal discrete control is $a^*(x(k)) = i, K^*(x(k)) = K_{i}(N-1)$ and can be defined by the following set of inequalities:
\begin{equation*}
\begin{cases}
x'(P_{i}(N-1) - P_{1}(N-1))x < 0\\
\vdots\\
x'(P_{i}(N-1) - P_{{i-1}}(N-1))x < 0\\
x'(P_{i}(N-1) - P_{{i+1}}(N-1))x \leq 0\\
\vdots\\
x'(P_{i}(N-1) - P_{{p}}(N-1))x \leq 0\\
\end{cases}
\end{equation*}
where the matrices are all symmetric. Note that $\Omega(k) = \{\Omega_i(N-1)\}_{i=1}^{p}$ is by definition a partition of $\mathbb{R}^n$. Each $\Omega_i(N-1)$ is associated with the discrete control $i$, the continuous state-feedback control $K_i(N-1)$ and the cost $P_i(N-1)$. Note that, at this step, each $\Omega_i(N-1)$ provides the optimal solution, i.e. $\iota(N-1) = \{1,\ldots,p\}$.


\textit{Step $N-2$:} Given the optimal solution $\Omega(N-1)$ of step $N-1$ we provide the optimal solution at step $N-2$. Let us write the cost-to-go function at step $N-2$:
\tiny
\begin{align*}
&J(x(N-2), u(N-2), a(N-2))=\\
&\min_{a(N-2),u(N-2)} E \Big\{x'(N-2)Mx(N-2)+u'(N-2)Ru(N-2)  + \\
&x'(N-1)P(N-1)x(N-1) \Big| \aleph_{N-2} \Big\},
\end{align*}
\normalsize
By linearity of the expected value we can write:
\tiny
\begin{align}\label{LQRsplit2}
&J(x(N-2),u(N-2),a(N-2))=\notag\\
&\min\limits_{a(N-2),u(N-2)} \Big\{ E \{x'(N-2)Mx(N-2) | \aleph_{N-2} \}+\notag\\
&E \{u'(N-2)Ru(N-2)   | \aleph_{N-2}\} +\notag\\
&E \{(x'(N-2)A'_{\sigma(N-2)}P(N-1)A_{\sigma(N-2)}x(N-2) | \aleph_{N-2}\} +\notag\\
&E \{(x'(N-2)A'_{\sigma(N-2)}P(N-1)B_{a(N-2)}u(N-2) | \aleph_{N-2}\} +\notag\\
&E \{(u'(N-2)B_{a(N-2)}'P(N-1)A_{\sigma(N-2)}x(N-2)  | \aleph_{N-2}\} +\notag\\
&E \{u'(N-2))B_{a(N-2)}'P(N-1)B_{a(N-2)}u(N-2) |\aleph_{N-2}\} \Big\}
\end{align}
\normalsize
Let us consider each addend of the rightside of \eqref{LQRsplit2}. Since $M$ and $R$ are constant matrices, $x(N-2)$ is $\aleph_{N-2}$-measurable, and $u(N-2)$ is not a random variable since it is the input that we choose to apply to the system at time $N-2$, the first two addends can be written as
\tiny
\begin{align*}
&E \{x'(N-2)Mx(N-2) +u'(N-2)Ru(N-2)| \aleph_{N-2}\}= \\
&x'(N-2)Mx(N-2) +u'(N-2)Ru(N-2).
\end{align*}
\normalsize
Moreover, $x'(N-2)Mx(N-2)$ does not depend on $a(N-2), u(N-2)$ and can be moved out of the min operator. Note that $P(N-1)$ depends, according to the optimal control policy $\Omega(N-1)$, on the random variable $x(N-1)$, which in turn depends on $\sigma(N-2), x(N-2), a(N-2), u(N-2)$. As a consequence, differently from step $N-1$, computing the expected values in Equation \eqref{LQRsplit2} is non-trivial. The main idea of this proof is to overcome this difficulty, instead of exploiting the conservative approximation used in \cite{VargasCDC2010}, by considering that such expected values can assume a finite number of values. In particular, given a discrete control input $a(N-2)=a$, a linear feedback $K(N-2)=K$ and a state $x(N-2)=x$ the probability that $P(N-1) = P_{i}(N-1)$ is equal to
\tiny
\begin{align*}
\sum\limits_{\sigma : (A_\sigma+B_a K)x \in \Omega_i(N-1)}\pi_\sigma.
\end{align*}
\normalsize
Let $\boldsymbol{\mu} \doteq \{\mu_{\sigma}\}_{\sigma \in \{1,\ldots, q\}}$ be a vector of $q$ natural numbers $\mu_{\sigma} \in \{1,\ldots, \omega(N-1) \}$ representing all possible combinations of sets $\Omega_i(N-1)$ towards which $x(N-2)$ can be driven by the occurrence of all $\sigma \in \Sigma$. It is easy to see that, given any $\boldsymbol{\mu}$, the expected values in Equation \eqref{LQRsplit2} can be computed. In particular, the third addend can be written as
\scriptsize
\begin{align}\label{eqAverage1}
&E \{(x'(N-2)A'_{\sigma(N-2)}P(N-1)A_{\sigma(N-2)}x(N-2) | \aleph_{N-2}\}=\notag\\
&x'(N-2) E \{A'_{\sigma(N-2)}P(N-1)A_{\sigma(N-2)}| \aleph_{N-2}\} x(N-2)=\notag\\
&x'(N-2) \left(\sum\limits_{\sigma=1}^{q} A'_\sigma P_{\mu_\sigma}(N-1) A_\sigma \pi_\sigma \right) x(N-2)
\end{align} 
\normalsize
Since $x(N-2)$ is $\aleph_{N-2}$-measurable, $P_i(N-1), i=1,\ldots,p$ are all symmetric and $\sigma(k)$ are i.i.d., the sum of the fourth and fifth addends can be written as
\tiny
\begin{align}\label{eqAverage2}
&2E \{(x'(N-2)A'_{\sigma(N-2)}P(N-1)B_{a(N-2)}u(N-2) | \aleph_{N-2}\}=\notag\\
&=2x'(N-2)E \{ A'_{\sigma(N-2)}P(N-1)| \aleph_{N-2}\}B_{a(N-2)}u(N-2)=\notag\\
&=2x'(N-2) \left(\sum\limits_{\sigma=1}^{q} A'_\sigma P_{\mu_\sigma}(N-1) \pi_\sigma \right) B_{a(N-2)}u(N-2).
\end{align}
\normalsize
The last addend can be written as:
\tiny
\begin{align}\label{eqAverage3}
&E \{u'(N-2))B_{a(N-2)}'P(N-1)B_{a(N-2)}u(N-2) |\aleph_{N-2}\}=\notag\\
&u'(N-2))B_{a(N-2)}' E \{ P(N-1) |\aleph_{N-2} \} B_{a(N-2)}u(N-2)=\notag\\
&u'(N-2))B_{a(N-2)}'  \left(\sum\limits_{\sigma=1}^{q} P_{\mu_\sigma}(N-1) \pi_\sigma \right) B_{a(N-2)}u(N-2).
\end{align}
\normalsize
We can now, for all states $x(N-2)$ that are driven by each $\sigma \in \Sigma$ to the set $\Omega_{\mu_\sigma}(N-1)$, rewrite \eqref{LQRsplit2} as follows:
\tiny
\begin{align} \label{costo_n-2}
&J(x(N-2),u(N-2),a(N-2), \boldsymbol\mu)=\notag\\
&\min\limits_{a(N-2)} \Big\{ x'(N-2)\left[M+ \left(\sum\limits_{\sigma=1}^{q} A'_\sigma P_{\mu_\sigma}(N-1) \pi_\sigma \right) \right]x(N-2)+\notag\\
&\min\limits_{u(N-2)} \Big\{ 2x'(N-2) \left(\sum\limits_{\sigma=1}^{q} A'_\sigma P_{\mu_\sigma}(N-1) \pi_\sigma \right) B_{a(N-2)}u(N-2)+\notag\\ 
& u'(N-2))B_{a(N-2)}' \left(\sum\limits_{\sigma=1}^{q} P_{\mu_\sigma}(N-1) \pi_\sigma \right) B_{a(N-2)}u(N-2)+\notag\\
&u'(N-2)Ru(N-2)\Big\}\Big\}
\end{align}
\normalsize
As a consequence the optimal linear feedback strategy can be computed as in step $N-1$ by:
\tiny
\begin{align}\label{eqOptimalKRecursiveStep}
&K_{a,\boldsymbol{\mu}}\doteq\notag\\
& -\left[R+ B_{a}' \left(\sum\limits_{\sigma=1}^{q} P_{\boldsymbol\mu_\sigma}(N-1) \pi_\sigma\right) B_{a} \right]^{-1} B_{a}' \left(\sum\limits_{\sigma=1}^{q} A'_{\sigma} P_{\boldsymbol\mu_\sigma}(N-1) \pi_\sigma\right) \bar A.
\end{align}
\normalsize
The set of states $x(N-2)$ such that $K_{a,\boldsymbol{\mu}}$ is indeed optimal is given by
\tiny
\begin{align*}
& \Gamma(a, \boldsymbol{\mu}) \doteq \{x \in \mathbb{R}^n : \forall \sigma \in \{1,\ldots,q\}, (A_\sigma + B_a K_{a,\boldsymbol{\mu}}) \in \Omega_{\mu_{\sigma}}(N-1)\}.
\end{align*}
\normalsize
Given the definition of $\Omega(N-1)$ of step $N-1$, then $\Gamma(a, \boldsymbol{\mu})$ can be defined by
\tiny
\begin{equation*}\label{eqGammaDef}
\begin{cases}
x'[(A_1 + B_a K_{a,\boldsymbol{\mu}})^{-1}]'(P_{\mu_{1}}(N-1) - P_{1}(N-1))(A_1 + B_a K_{a,\boldsymbol{\mu}})^{-1}x < 0\\
\vdots\\
x'[(A_1 + B_a K_{a,\boldsymbol{\mu}})^{-1}]'(P_{\mu_{1}}(N-1) - P_{{p}}(N-1))(A_1 + B_a K_{a,\boldsymbol{\mu}})^{-1}x \leq 0\\
\vdots\\
x'[(A_q + B_a K_{a,\boldsymbol{\mu}})^{-1}]'(P_{\mu_{q}}(N-1) - P_{1}(N-1))(A_q + B_a K_{a,\boldsymbol{\mu}})^{-1}x < 0\\
\vdots\\
x'[(A_q + B_a K_{a,\boldsymbol{\mu}})^{-1}]'(P_{\mu_{q}}(N-1) - P_{{p}}(N-1))(A_q + B_a K_{a,\boldsymbol{\mu}})^{-1}x \leq 0
\end{cases}
\end{equation*}
\normalsize
Note that, for any given $a \in A$, $\Gamma(a, \boldsymbol{\mu})$ is a subset of $\mathbb{R}^n$: its complement $\Gamma^C(a, \boldsymbol{\mu})$ can be easily defined by replacing in \eqref{eqGammaDef} $<$ and $\leq$ respectively with $\geq$ and $>$, and represents a set where the optimal linear feedback strategy cannot be computed as in \eqref{eqOptimalKRecursiveStep}. Also, for any given $a \in A$ and any $\boldsymbol{\mu}, \boldsymbol{\mu'}$, the intersection $\Gamma(a, \boldsymbol{\mu}) \cap \Gamma(a, \boldsymbol{\mu'})$ is not necessarily empty. Define now the set $\{\Omega_i\}_{i=1}^{\omega}$ of disjoint sets given by all possible combinations of intersections of sets $\Gamma(a,\boldsymbol{\mu}), a \in A, \boldsymbol{\mu}$, and such that for all $\bar a \in A$ at least a set $\Gamma(\bar a,\boldsymbol{\bar\mu})$ belongs to the intersection. As a consequence each $\Omega_i$ can be defined as a finite intersection of sets $\Gamma(a, \boldsymbol{\mu})$ and $\Gamma^C(a', \boldsymbol{\mu'})$, and can be therefore characterised by a finite set of quadratic inequalities. Consider now a generic set
\tiny
$$
\Omega_i = \Gamma(a_1,\boldsymbol{\mu_1}) \cap \Gamma(a_2,\boldsymbol{\mu_2}) \cap \cdots \cap \Gamma(a_z,\boldsymbol{\mu_z}) \cap \Gamma^C(a_{z+1},\boldsymbol{\mu_{z+1}}) \cap \cdots \cap \Gamma^C(a_\gamma,\boldsymbol{\mu_\gamma}).
$$
\normalsize
For any $x(N-2) \in \Omega_i$ we have the choice to apply a finite number of optimal control pairs $\{(a_j, K_{a_j, \mu_j})\}_{ j \in \{1,\ldots,z\}}$, each associated to the optimal cost
\tiny
\begin{align*}
& J(x(N-2), a_j, \boldsymbol\mu_j) \doteq x'(N-2) P_{a_j, \boldsymbol\mu_j}(N-2) x(N-2).
\end{align*}
\normalsize 
To choose the optimal control feedback within $\Omega_i$ we partition it, as in step $N-1$, with a finite collection of disjoint sets $\{\Omega_{i,a,\boldsymbol\mu}\}$ defined by
\tiny
\begin{align*}
\Omega_{i,a,\boldsymbol\mu} \doteq \{x(N-2) \in \Omega_i : a, \boldsymbol\mu = \arg\min\limits_{a,\boldsymbol\mu} J(x(N-1), a, \boldsymbol\mu) \}.
\end{align*}
\normalsize
$\Omega_{i,a,\boldsymbol\mu}$ is the set of all states $x(N-2)$ such that the corresponding optimal control is $a^*(x(k)) = a, K^*(x(k)) = K_{a,\boldsymbol\mu}$. As a consequence $\Omega_{i,a_j,\boldsymbol\mu_j}$ can be defined as follows:
\begin{equation}\label{eqPartitionRecursiveStep}
\Omega_{i} \cap
\begin{cases}
x'(P_{a_j,\boldsymbol\mu_j}(N-1) - P_{a_1,\boldsymbol\mu_1}(N-1))x < 0\\
\vdots\\
x'(P_{a_j,\boldsymbol\mu_j}(N-1) - P_{a_{j-1},\boldsymbol\mu_{j-1}}(N-1))x < 0\\
x'(P_{a_j,\boldsymbol\mu_j}(N-1) - P_{a_{j+1},\boldsymbol\mu_{j+1}}(N-1))x \leq 0\\
\vdots\\
x'(P_{a_j,\boldsymbol\mu_j}(N-1) - P_{a_z,\boldsymbol\mu_z}(N-1))x \leq 0\\
\end{cases},
\end{equation}
where the matrices are all symmetric. Note that $\{\Omega_{i,a_j,\boldsymbol\mu_j}\}_{j=1}^{z}$ is by definition a partition of $\Omega_i$. To each $\Omega_{i,a_j,\boldsymbol\mu_j}$ is associated the control $a_j$, $K_{a_j,\boldsymbol\mu_j}$ and the cost $P_{a_j,\boldsymbol\mu_j}$. Applying the same partition to each $\Omega_i$ provides the optimal solution for any state in the set $\Omega \doteq \bigcup\limits_{i,a,\boldsymbol\mu} \Omega_{i,a,\boldsymbol\mu} \subseteq \mathbb{R}^n$. The set of states $\Omega^C \doteq \mathbb{R}^n \setminus \Omega$ such that the optimal solution cannot be computed using the method above can be also defined as a finite union of sets of quadratic inequalities. To compute a suboptimal solution for states $x(N-2) \in \Omega^C(N-2)$ we can partition it in a finite number of sets $\Omega_{a,\boldsymbol\mu}^C$ as in \eqref{eqPartitionRecursiveStep} according to the set of all pairs $a, K_{a, \boldsymbol\mu}$. Of course this will not be the optimal solution, since by definition there exists at least a discrete control $a$ such that $K_{a, \mu}$ is not the corresponding optimal linear feedback. We can define $\Omega(N-2)$ by the union of all sets that partition $\Omega$ and $\Omega^C$, as defined above. $\Omega(N-2)$, together with the associated discrete controls, linear feedback gains and costs, will be provided as input to step $N-3$. Also, the set $\iota(N-2)$ can be easily characterised since only the solutions for the sets partitioning $\Omega$ are optimal. The solution at any other step $k = 1, \ldots, N-3$ given the solution at step $k+1$ can be obtained iterating the same reasoning of step $N-2$. The only difference is that a set $\Omega_{i, a, \boldsymbol\mu}$ is associated to an optimal solution only if $\Omega_{\mu_{\sigma}}(k+1)$ is associated to an optimal solution for all $\sigma \in \Sigma$. Once arrived to the initial step $N=0$, we can define the sets of initial conditions such that the strategy is optimal. This concludes the proof. $\qed$
\end{theorem}


\section{Conclusion}

This paper proposes a novel paradigm of networked control where the exploitation of redundancy in routing actuation data tremendously improves the control performance, by exploiting in an optimal way the advantages of different paths characterised by propagation characteristics that are at odds one another. We show that the co-design of controller and routing redundancy, given a routing that is defined a priori, strongly improves the control performance. We also provide a methodology to co-design controller and dynamic routing redundancy: note that, according to our solution, routing depends on the plant state, and it is even possible to decide not to route any actuation data for some time instants, which generates an \emph{event-triggered} control strategy to send actuation data not at all time steps, but only when necessary according to the current plant state. In future work we plan to implement a tool that applies our algorithms and to extend our methodology to more general classes of stochastic systems. Also, we will apply the methodologies in this paper to address the dynamic scheduling and routing co-design in communication protocols for wireless control systems such as WirelessHART and ISA100 (see \cite{AlurTAC11} and \cite{DInnocenzoTAC13} for details).




\bibliography{mcnbib}             

\begin{thebibliography}{23}
\providecommand{\natexlab}[1]{#1}
\providecommand{\url}[1]{\texttt{#1}}
\providecommand{\urlprefix}{URL }
\expandafter\ifx\csname urlstyle\endcsname\relax
  \providecommand{\doi}[1]{doi:\discretionary{}{}{}#1}\else
  \providecommand{\doi}{doi:\discretionary{}{}{}\begingroup
  \urlstyle{rm}\Url}\fi

\bibitem[{Astr\"om and Wittenmark(1997)}]{Astrom97j1}
Astr\"om, K. and Wittenmark, B. (1997).
\newblock \emph{{Computer-Controlled Systems: Theory and Design}}.
\newblock Prentice Hall.

\bibitem[{Bertsekas(2005)}]{BertsekasOptStocCon}
Bertsekas, D.P. (2005).
\newblock \emph{Dynamic programming and optimal Control Vol. I \& II}.
\newblock Athena Scientific, Belmont, Massachussets.

\bibitem[{Borrelli et~al.(2005)Borrelli, Baotic, Bemporad, and
  Morari}]{BemporadAutomatica2005}
Borrelli, F., Baotic, M., Bemporad, A., and Morari, M. (2005).
\newblock Dynamic programming for constrained optimal control of discrete-time
  linear hybrid systems.
\newblock \emph{Automatica}, 41(4), 1709--1721.

\bibitem[{Costa et~al.(2005)Costa, Marques, and
  Fragoso}]{costa_discrete-time_2005}
Costa, O.L.V., Marques, R.P., and Fragoso, M.D. (2005).
\newblock \emph{{Discrete-Time} Markov Jump Linear Systems}.
\newblock Springer.

\bibitem[{D'Innocenzo et~al.(2013)D'Innocenzo, {M.D. Di Benedetto}, and
  Serra}]{DInnocenzoTAC13}
D'Innocenzo, A., {M.D. Di Benedetto}, and Serra, E. (2013).
\newblock Fault tolerant control of multi-hop control networks.
\newblock \emph{IEEE Transactions on Automatic Control}, 58(6), 1377--1389.

\bibitem[{Gatsis et~al.(2014)Gatsis, A.Ribeiro, and Pappas}]{PappasTAC2014}
Gatsis, K., A.Ribeiro, and Pappas, G.J. (2014).
\newblock Optimal power management in wireless control systems.
\newblock \emph{IEEE Transactions on Automatic Control}, 59(6), 1495--1510.

\bibitem[{{G.C. Walsh} et~al.(2002){G.C. Walsh}, Ye, and {L.G.
  Bushnell}}]{walsh_stability_2002}
{G.C. Walsh}, Ye, H., and {L.G. Bushnell} (2002).
\newblock {Stability Analysis of Networked Control Systems}.
\newblock \emph{{IEEE} Transactions on Control Systems Technology}, 10(3),
  438--446.
\newblock \doi{10.1109/87.998034}.

\bibitem[{Gupta et~al.(2009)Gupta, Dana, Hespanha, Murray, and
  Hassibi}]{MurrayTAC2009}
Gupta, V., Dana, A., Hespanha, J., Murray, R., and Hassibi, B. (2009).
\newblock Data transmission over networks for estimation and control.
\newblock \emph{{IEEE} Transactions on Automatic Control}, 54(8), 1807--1819.

\bibitem[{Han et~al.(2010)Han, Xiuming, Aloysius, Nixon, Blevins, and
  Chen}]{SongIECON2010}
Han, S., Xiuming, Z., Aloysius, K., Nixon, M., Blevins, T., and Chen, D.
  (2010).
\newblock Control over wirelesshart network.
\newblock In \emph{IECON 2010 - 36th Annual Conference on IEEE Industrial
  Electronics Society}, 2114--2119.

\bibitem[{{I.F. Akyildiz} and {I.H. Kasimoglu}(2004)}]{akyildiz_wireless_2004}
{I.F. Akyildiz} and {I.H. Kasimoglu} (2004).
\newblock {Wireless Sensor and Actor Networks: Research Challenges}.
\newblock \emph{{Ad Hoc Networks}}, 2(4), 351--367.

\bibitem[{{J.P. Hespanha} et~al.(2007){J.P. Hespanha}, Naghshtabrizi, and
  Xu}]{Hespanha2007}
{J.P. Hespanha}, Naghshtabrizi, P., and Xu, Y. (2007).
\newblock {A Survey of Recent Results in Networked Control Systems}.
\newblock \emph{Proceedings of the IEEE}, 95(1), 138--162.

\bibitem[{{K.-E. {\AA}rz\'{e}n} et~al.(2006){K.-E. {\AA}rz\'{e}n}, {A. Bicchi},
  {S. Hailes}, {K. H. Johansson}, and {J. Lygeros}}]{Arzen06}
{K.-E. {\AA}rz\'{e}n}, {A. Bicchi}, {S. Hailes}, {K. H. Johansson}, and {J.
  Lygeros} (2006).
\newblock On the design and control of wireless networked embedded systems.
\newblock In \emph{Proceedings of the 2006 IEEE Conference on Computer Aided
  Control Systems Design, Munich, Germany}, 440--445.

\bibitem[{Lancaster and Rodman(1995)}]{LancasterRiccati}
Lancaster, P. and Rodman, L. (1995).
\newblock \emph{Algebraic Riccati Equations}.
\newblock Clarendon Press - Oxford.

\bibitem[{Matei et~al.(2008)Matei, Martins, and Baras}]{BarasProcIFAC}
Matei, I., Martins, N., and Baras, J. (2008).
\newblock Optimal linear quadratic regulator for markovian jump linear systems,
  in the presence of one time-step delayed mode observations.
\newblock In \emph{Proceedings of the 17 World Congress The International
  Federation of Automatic Control}.

\bibitem[{{M.C.F. Donkers} et~al.(2011){M.C.F. Donkers}, {W.P.M.H. Heemels},
  {Nathan van de Wouw}, and {Laurentiu Hetel}}]{HeemelsTAC11}
{M.C.F. Donkers}, {W.P.M.H. Heemels}, {Nathan van de Wouw}, and {Laurentiu
  Hetel} (2011).
\newblock {Stability Analysis of Networked Control Systems Using a Switched
  Linear Systems Approach}.
\newblock \emph{IEEE Transactions on Automatic Control}, 56(9), 2101 --2115.
\newblock \doi{10.1109/TAC.2011.2107631}.

\bibitem[{Mesquita et~al.(2012)Mesquita, Hespanha, and Nair}]{MHNAUT12}
Mesquita, A., Hespanha, J., and Nair, G. (2012).
\newblock Redundant data transmission in control/estimation over lossy
  networks.
\newblock \emph{Automatica}, 48, 1020--1027.

\bibitem[{Pajic et~al.(2011)Pajic, Sundaram, Pappas, and
  Mangharam}]{PajicTAC2011}
Pajic, M., Sundaram, S., Pappas, G., and Mangharam, R. (2011).
\newblock The wireless control network: a new approach for control over
  networks.
\newblock \emph{IEEE Transactions on Automatic Control}, 56(10), 2305--2318.

\bibitem[{{R. Alur} et~al.(2011){R. Alur}, {A. D'Innocenzo}, {K.H. Johansson},
  {G.J. Pappas}, and {G. Weiss}}]{AlurTAC11}
{R. Alur}, {A. D'Innocenzo}, {K.H. Johansson}, {G.J. Pappas}, and {G. Weiss}
  (2011).
\newblock {Compositional Modeling and Analysis of Multi-Hop Control Networks}.
\newblock \emph{IEEE Transactions on Automatic Control, Special Issue on
  Wireless Sensor and Actuator Networks}, 56(10), 2345--2357.

\bibitem[{Smarra et~al.(2015)Smarra, D'Innocenzo, and {Di Benedetto,
  M.D.}}]{SmarraECC15}
Smarra, F., D'Innocenzo, A., and {Di Benedetto, M.D.} (2015).
\newblock Approximation methods for optimal network coding in a multi-hop
  control network with packet losses.
\newblock In \emph{Proceedings of the $14^{th}$ European Control Conference
  (ECC'15), Linz, Austria, July 15-17}.

\bibitem[{Tabbara et~al.(2007)Tabbara, Ne\v{s}i\'{c}, and
  Teel}]{TabbaraTAC2007}
Tabbara, M., Ne\v{s}i\'{c}, D., and Teel, A. (2007).
\newblock {Stability of Wireless and Wireline Networked Control Systems}.
\newblock \emph{IEEE Transactions on Automatic Control}, 52(7), 1615--1630.

\bibitem[{Vargas et~al.(2010)Vargas, Ishihara, and do~Val}]{VargasCDC2010}
Vargas, A.N., Ishihara, J.Y., and do~Val, J. (2010).
\newblock Linear quadratic regulator for a class of markovian jump systems with
  control in jumps.
\newblock \emph{Conference on Decision and Control}.

\bibitem[{{W. Zhang} et~al.(2001){W. Zhang}, {M.S. Branicky}, and {S.M.
  Phillips}}]{Zhang2001}
{W. Zhang}, {M.S. Branicky}, and {S.M. Phillips} (2001).
\newblock {Stability of Networked Control Systems}.
\newblock \emph{IEEE Control Systems Magazine}, 21(1), 84--99.
\newblock \doi{10.1109/37.898794}.

\bibitem[{Zhang et~al.(2009)Zhang, Abate, and Hu}]{AbateACC2009}
Zhang, W., Abate, A., and Hu, J. (2009).
\newblock Efficient suboptimal solutions of switched lqr problems.
\newblock \emph{American Control Conference}, 37(4), 1084 -- 1091.

\end{thebibliography}
                                               

%
%
%
%
%

\end{document}